\documentstyle[12pt,a4]{article}  
  
\def\Box{\vbox{\hrule%  
\hbox{\vrule height 5pt depth 0pt\kern 5pt\vrule height 5pt depth   
0pt}%  
\hrule}}  
  
\newtheorem{theo}{Theorem}    
\newtheorem{prop}[theo]{Proposition}  
\newtheorem{coro}[theo]{Corollary}  
\newtheorem{lemme}[theo]{Lemma}

\def\NN{\hbox{\sf I\hskip -1pt N}}  
\def\pNN{\hbox{{\scriptsize \sf I}\hskip -1pt {\scriptsize \sf N}}}  
  
\def\RR{\hbox{\sf I\hskip -2pt R}}  
\def\pRR{\hbox{{\scriptsize \sf I}\hskip -1pt {\scriptsize \sf R}}}  
\def\TT{\hbox{\sf T\hskip -4pt T}}  
\def\pTT{\hbox{{\scriptsize \sf T}\hskip -3pt {\scriptsize \sf T}}}  
\def\ZZ{\hbox{\sf Z\hskip -4pt Z}}  
\def\pZZ{\hbox{{\scriptsize \sf Z}\hskip -4pt {\scriptsize \sf Z}}}  
\def\UU{{\bf U}}        
\def\CC{{\bf C}}  
  
\def\card{\#}

\def\sup{\hbox{\rm Sup}}  
\def\mvide{{\epsilon }}  
\def\mod{\hbox{\rm mod}}

\def\d{{\rm d}}  
\def\A{{\cal A}}  
\def\B{{\cal B}}

\def\H{{\cal H}}  
\def\L{{\cal L}}  
\def\N{{\cal N}}  
\def\R{{\cal R}}

\def\x{{\bf\sf x}}

\begin{document}  
  
\title  
{  
\bf  Ergodic averages with deterministic weights  
}  
  
\author{  
Fabien Durand, Dominique Schneider}  

\date{}  

\maketitle  
  
\vskip 1.5cm  

Math. classification: 37A05, 28D05, 11K99.

\medskip


\section{Introduction.}  
  
The purpose of this paper is to study ergodic averages with deterministic
weights. More precisely we study the convergence
of the ergodic averages of the type
$\frac{1}{N} \sum_{k=0}^{N-1} \theta (k) f \circ T^{u_k}$ where $\theta
= (\theta (k) ; k\in \NN)$ is a bounded sequence and $u = (u_k ; k\in \NN)$
a strictly increasing sequence of integers such that for some $\delta<1$
$$
S_N (\theta , u) := \sup_{\alpha \in
\pRR} \left|   \sum_{k=0}^{N-1} \theta (k) \exp (2i\pi \alpha u_k )
\right| = O (N^{\delta}) \ ,
\leqno{({\cal H}_1)}  
$$ 
i.e., there exists a
constant $C$ such that $S_N (\theta , u) \leq C N^{\delta} $. We define
$\delta (\theta , u)$ to be the infimum of the $\delta $  satisfying
$\H_1$ for  $\theta $  and $u$.

About ${\cal H}_1$, in the
case where $\theta $ takes its values in $\UU$ (the set of complex numbers of modulus
1), it is clear that for all sequences $\theta$ and $u$, $\delta (\theta
, u )$ is smaller than or equal to 1 and it is well-known (see \cite{Ka}
for example) that it is greater than or
equal to $1/2$. Few explicit sequences $\theta$ are known to have
$\delta (\theta , u)$ strictly smaller than 1. 

When $u_k = k$, for all $k\in \NN$, we know \cite{Ru,Sh} 
that for the Rudin-Shapiro sequence (and its generalizations \cite{AL,MT}) we have
$\delta (\theta , u) = 1/2$.  For the Thue-Morse sequence $\delta (\theta ,
u) =  (\log 3 )/ (\log 4)$ \cite{G}. When $\theta$ is a
$q$-multiplicative sequence we will give a way to construct sequences
fulfilling ${\cal H}_1$.

When $u$ is a subsequence of $\NN$ we will also
give some examples of sequences $\theta$ satisfying ${\cal H}_1$. More
attention will be payed to the special case $u_k = k + v_k$, $k\in
\NN$, where $v= (v_k ; k\in \NN )$ is non-decreasing with $v_k = O
(k^{\varepsilon})$, $\varepsilon <1$.

We say $(X,\B , \mu , T)$ is a {\it dynamical system} if $(X,\B , \mu)$
is a probabilistic space, $T$ is a measurable map from $X$ to $X$ and
$\mu$ is $T$-invariant. A {\it good sequence for the pointwise ergodic
theorem in} $L^p (\mu )$, $p\geq 1$, is an increasing integer sequence $(u_n ; n\in \NN)$ such that
for all dynamical systems $(X,\B , \mu , T)$, for all $f\in L^p (\mu )$,
we have

$$  
\mu \left\{ x\in X ; \lim_{N\rightarrow +\infty } \frac{1}{N}   
\sum_{k=0}^{N-1}  f (T^{u_k} x)  \hbox{ exists } \right\} = 1 .  
$$  

Our main result is the following.

\begin{theo}  
\label{thergo}  
Let $\theta = (\theta (n) ; n\in \NN )$ be a bounded sequence of complex numbers and   
$u = (u_n ; n\in \NN )$ be a strictly increasing sequence of
 integers. Suppose that Condition ${\H_1}$ is satisfied.

Then, for any dynamical system $(X,\B , \mu , T)$ and any $f\in L^2   
(\mu )$ we have   
$$  
\mu \left\{ x\in X ; \lim_{N\rightarrow +\infty } \frac{1}{N}   
\sum_{k=0}^{N-1} \theta (k) f (T^{u_k} x)  = 0 \right\} = 1 .  
$$  
Moreover, for all $\varepsilon > 0$, all $f\in L^{2+\varepsilon} (\mu )   
$ and for all $\beta > (\delta + 2)/3$ we have   
$$  
\mu \left\{ x\in X ; \lim_{N\rightarrow +\infty } \frac{1}{N^{\beta}}   
\sum_{k=0}^{N-1} \theta (k) f (T^{u_k} x)  = 0 \right\} = 1 .  
$$  
\end{theo}  
Remark that if $u$ is a good sequence for the pointwise ergodic theorem in $L^1(\mu)$   
then the first conclusion is satisfied for all $f \in L^1 (\mu)$.

Under the stronger Condition $\H_2$ given below we can give more
 information about the speed
 of convergence when $f$ belongs to $L^2 (\mu )$. Let $\theta = (\theta
 (k) ; k\in \NN)$ and $u = (u_k ; k\in \NN)$ be as in Theorem
 \ref{thergo}. We say they satisfy Condition $\H_2$ if for some $0 < \rho < 1$ and $D>0$ we have
$$
\sup_{\alpha \in \pRR} 
\left|   
\sum_{k=m}^{n} \theta (k) \exp (2i\pi \alpha u_k )
\right| 
\leq
D(n-m)^\rho
\zeta (n) ,
\leqno{( \H_2)}   
$$
for all $m\leq n$ where $\zeta (n) = O (n^\epsilon)$ for all $\epsilon >
0$. We remark that $\H_2$ implies $\H_1$. When $\theta$ is
a $q$-multiplicative sequence we prove that Conditions $\H_1$
and $\H_2$ are equivalent.

We will show (Section \ref{cororot}) in the special case of a rotation dynamical system we have
some more precise results than those of Theorem \ref{thergo}. We will
exhibit a large class of functions for which the speed of convergence is
uniform in $x$. In this case and when $\theta (k) = 1$ and $u_k = k$ (i.e., the standard   
ergodic mean) we knew, using the Denjoy-Koksma inequality, that we can
obtain the speed of convergence whenever $f$ is has bounded variations
(see \cite{KN} or \cite{GS}).  

For $u = (k^2 ; k\in \NN )$ we do not know sequences that satisfy
${\cal H}_1$. But when $\theta$ is a $q$-multiplicative sequence we prove
the following result.

\begin{prop}  
\label{carre}
Let $\theta $ be a $q$-multiplicative sequence with empty
 spectrum. Then, for all real number $\alpha $  
$$  
\lim_{N\rightarrow +\infty}  
\frac{1}{N} \sum_{k=1}^{N} \theta (k) e ( {k^2} \alpha ) = 0 .  
$$  
\end{prop}  

The main difficulty in the proof of this result is when $\alpha $ is an
irrational. In this case we use an ergodic approach and the van der
Corput inequality.

For more information about Proposition \ref{carre} we invite the
interested reader to look at the article of M. Mend\`es France \cite{M},
especially Corollaire 2.

As noticed by the referee Proposition \ref{carre} can be extended to the
case where $k^2$ is replaced by any polynomial of degree $d\geq 2$ with
rational coefficients. We can prove this result using an induction on the
degree of the polynomial and the van der Corput inequality. Then it can
be easily extended to
the case where $\alpha k^2$ is replaced by any polynomial of degree
$d\geq 2$ with at least one irrational coefficient.

The following theorem investigates the statistics of dynamical
systems. We prove a central limit theorem for some weighted ergodic
means in the case of the rotations. First we need two definitions. 

Let $(Z_n ; n\in \NN )$ be a sequence of real ramdom variables defined
on the probability space $(X, \B , \mu )$. We say that $(Z_n ; n\in \NN
)$ converges in law to the Gaussian random variable $\N (0,1)$, and we write it $
Z_n   \Longrightarrow^{\hskip -12pt \L}_{\hskip +7pt } \N (0,1)$, if
$$
\lim_{n\rightarrow +\infty } \mu \{ Z_n \geq t \} = \frac{1}{\sqrt{2\pi}}\int_{t}^{+\infty}
e^{-\frac{u^2}{2}} du .
\leqno{\forall t\in \RR ,}  
$$  

We define ${\rm lip} (a)$ to be the set of functions $f\in L_2 (\mu)$ such that
$| f(x) - f(y) | \leq |x-y|^a $.

\begin{theo}  
\label{cvlaw}  
Let $(\theta (n) ; n\in \NN ) \in \{ -1 , 1 \}^{\pNN}$ and $(u_n ; n\in \NN )$ satisfying  Condition ${\cal H}_1$. Let $\theta^{+}$ and $\theta^{-}$ be the sequences in $\{ 0   
, 1\}^{\pNN}$ defined by $\theta^{+} (k) = (\theta (k) +1 )/2 $ and   
$\theta^{-} (k) = (1 - \theta (k))/2 $. Let $(\TT,\B,\mu , R_{\alpha })$   
be a rotation dynamical system where $\alpha $ is an irrational number.   
For all $\beta \in ] \delta , 1[ $, if there exist $\sigma > 1-\beta$ and $\zeta \in ]0,1]$ such that $|\frac{u_n}{N} - \zeta | = O (N^{-\sigma})$ then there exists a continuous function   
$\tilde{f}$ on $\TT$ such that    
$$  
\frac{1}{N^{\beta}} \sum_{k=0}^{N-1} \theta^{+} (k) \tilde{f} \circ R_{\alpha}^{u_k}   \Longrightarrow^{\hskip -12pt \L}_{\hskip +7pt } \N (0,1)   
.    
$$  
Moreover, if $d$ is the Diophantine type of $\alpha$ then $\tilde{f}$   
is ${\rm lip} (a)$ with $a< (1-\beta)/d$.  
  
We have the same conclusions for $\theta^{-}$.  
\end{theo}

This work is divided into four sections. In Section 2 we prove Theorem \ref{thergo}
and Theorem \ref{cvlaw} and we deduce some corollaries. For example we
remark that the conclusions of Theorem \ref{thergo} also hold when
$(\theta(k) ; k\in \NN )$ is a non-bounded centered sequence of i.i.d. random
variables with a finite second moment (see \cite{SW}). We also make some comments in the case there is
no dynamical system structure. The third section is devoted to the
$q$-multiplicative sequences. We recall some results established in
\cite{LMM} and we give an efficient sufficient condition
for $q$-multiplicative sequences taking values in a finite set to fulfill
condition ${\cal H}_1$. In the last section we prove Proposition
\ref{carre} and we obtain further results about ergodic averages in the case where $\theta$ is a $q$-multiplicative sequence.

\section{Convergence of ergodic weighted averages under Condition ${\cal H}_1$}  
  
In what follows we will write $e(x)$ instead of $\exp (2i\pi x)$.  
  
\subsection{Proof of Theorem \ref{thergo}}  
  
We start with the proof of the second conclusion. We first prove it when $N$ tends to infinity along some subset of $\NN$. Then we come back to $\NN$.  
  
\medskip  
  
Let $\varepsilon > 0 $ and  $f\in L^{2+\varepsilon   
} (\mu )$. From the Spectral Lemma (see \cite{Kr}) we have  
$$  
\left| \left|   
\frac{1}{N^{\beta}} \sum_{k=0}^{N-1} \theta (k) f \circ T^{u_k}   
\right| \right|_{2,\mu}  
=  
\left| \left|  
\frac{1}{N^{\beta}} \sum_{k=0}^{N-1} \theta (k) e( \alpha u_k )  
\right| \right|_{2, \pTT,\mu_f ({\rm d}\alpha )}   
\leq   
\frac{C}{N^{\beta - \delta}} \left| \left|  
f  
\right| \right|_{2,\mu} ,  
$$  
where $\mu_f$ is the spectral measure of $T$ at $f$. We set $\beta_0 =   
(\delta +2 )/ 3$ and we take $\beta > \beta_0$. Hence we can choose $\gamma \in \RR$ such that $1/(\beta   
- \delta) < \gamma < 1/2(1-\beta )$. We set $\N_{\gamma}  = \{ [ n^{\gamma} ] ; n\geq 1   
\}$, where $[.]$ is the integer part function. From the choice of $\gamma $ it comes that  
$$  
\sum_{N\in \N_{\gamma}}  
\left| \left|   
\frac{1}{N^{\beta}} \sum_{k=0}^{N-1} \theta (k) f \circ T^{u_k}   
\right| \right|_{2,\mu}  
< \infty .  
$$  
Consequently  
$$  
\mu \left\{ x\in X ; \lim_{N\rightarrow +\infty , N\in \N_{\gamma} }   
\frac{1}{N^{\beta}} \sum_{k=0}^{N-1} \theta (k) f (T^{u_k} x)  = 0 \right\}   
= 1 .  
$$  
  
Now we come back to the whole set $\NN$. There exists $K$ such that
$|\theta (k)|<K$ for all $k\in \NN$. Let $M > 0$, sufficiently large, and let $N$ be the unique integer such that $[N^{\gamma}] \leq M <   
[(N+1)^{\gamma}$]. Then  
$$  
\left|   
\frac{1}{M^{\beta}} \sum_{k=0}^{M-1} \theta (k) f (T^{u_k} x)   
\right|  
\leq  
\left|  
\frac{2}{N^{\beta \gamma}} \sum_{0\leq k < N^{\gamma }} \theta (k) f   
(T^{u_k} x)   
\right|  
+  
\frac{K}{N^{\beta \gamma }} \sum_{N^{\gamma } \leq k <(N+1)^{\gamma }}    
\left| f (T^{u_k} x) \right| .  
$$  
  
From the first step of the proof, the first term tends to zero. We have   
to prove the second term also tends to $0$. Without loss of generality   
we can suppose $f$ is positive. We set $f = f_1 + f_2$ where $f_1 = f .  
\chi_{\{ f\leq \sqrt{N} \} }$ and $f_2 = f . \chi_{ \{ f > \sqrt{N} \} }$, $\chi_A$ being the characteristic function of the set $A$. We have   
$$  
\frac{1}{N^{\beta \gamma}} \sum_{N^{\gamma } \leq k < (N+1)^{\gamma }}   
 f_1 (T^{u_k} x)   
\leq  
\frac{1}{N^{\beta \gamma}} \sum_{N^{\gamma } \leq k < (N+1)^{\gamma }}   
\sqrt{N}  
$$  
$$  
\leq  
B \frac{N^{\gamma - 1}}{N^{\beta \gamma }} \sqrt{N}  
=  
B N^{(1-\beta) \gamma -\frac{1}{2}} \longrightarrow_{N\rightarrow +\infty } 0,  
$$  
where $B$ is such that $(N+1)^{\gamma} - N^{\gamma} \leq B N^{\gamma -
1}$, for all $N\in \NN$. Now we look what happens with $f_2$. To prove that   
$$
\mu \left\{ x\in X ; \lim_{N\rightarrow +\infty }   
\frac{1}{N^{\beta \gamma}} \sum_{
N^{\gamma } \leq k < (N+1)^{\gamma }
} f_2 (T^{u_k} x)  = 0 \right\}   
= 1   
$$ 
it suffices to prove that $\int_X \frac{1}{N^{\beta \gamma}} \sum_{N^{\gamma } \leq k < (N+1)^{\gamma }} f_2 (T^{u_k} x) d\mu$ is the general term of a convergent series. The measure $\mu$ being $T$-invariant we have  
$$  
\int_X   
\frac{1}{N^{\beta \gamma }} \sum_{N^{\gamma } \leq k < (N+1)^{\gamma }}    
f_2 (T^{u_k} x) \d\mu  
=
\frac{1}{N^{\beta \gamma }} \sum_{N^{\gamma } \leq k < (N+1)^{\gamma }}  
\int_X  
f_2 (T^{u_k} x) \d\mu  
$$  
$$  
\leq  
\frac{1}{N^{\beta \gamma }}   
B N^{\gamma - 1}  
\int_X  
f_2 ( x) \d\mu  
=  
B  
\int_X  
\frac{f_2 ( x)}{N^{1 - \gamma ( 1 - \beta )}} \d\mu  
$$  
We set $\alpha = 1-\gamma (1-\beta )$ and  we will study the series  
$$  
S=  
\sum_{N\geq 1}  
\int_X  
\frac{f_2 ( x)}{N^{\alpha}} \d\mu .  
$$  
Let $p = (1+\varepsilon / 2)/\alpha $ and $q = (1 - 1/p)^{-1}$, the H\"older inequality gives  
$$  
S \leq  
\left(  
\int_X  
\sum_{N\geq 1}  
\left(  
\frac{f ( x)}{N^{\alpha}}  
\right)^{p}  
\d\mu  
\right)^{\frac{1}{p}}  
\left(  
\sum_{N\geq 1}  
\mu \{ f^2 > N \}  
\right)^{\frac{1}{q}}  
$$  
$$  
=  
\left(  
\int_X  
|f ( x)|  
^{p}  
\d\mu  
\right)^{\frac{1}{p}}  
\left(  
\sum_{N\geq 1}  
\frac{1}{N^{p\alpha}}  
\right)^{\frac{1}{p}}  
\left(  
\sum_{N\geq 1}  
\mu \{ f^2 > N \}  
\right)^{\frac{1}{q}} .  
$$  
The last series converges because $f$ belongs to $L^2 (\mu )$. The second also converges because $\alpha p = 1+\varepsilon / 2 $. Moreover   
$$  
p   
=   
\frac{1+\frac{\varepsilon}{2}}{\frac{1}{2} + \frac{1}{2} -\gamma (1-\beta )}  
\leq  
2+ \varepsilon ,  
$$  
consequently $S$ is finite. This ends the proof of the second conclusion.  
  
\medskip  
  
The same computation with $f\in L^2 (\mu )$, $\beta = 1$, $\gamma >1/(1-\delta)$ and $p=q=2$ allows us to obtain the first conclusion.  
\hfill $\Box$

\medskip  
  
{\bf Remark.} The first conclusion of Theorem \ref{thergo} can be extended to every $f\in L^p (\mu ) $ with $p>1$. On the other hand if $(u_k ; k\in \NN)$ is a good sequence for the Ergodic Theorem in $L^1 (\mu )$ then we know from Banach Principle (see \cite{Kr}) that the set of functions for which we have almost sure convergence (this set contains $L^2 (\mu)$) is closed in $L^1 (\mu )$. Hence, $L^2 (\mu )$ being dense in $L^1 (\mu )$, the first conclusion holds in $L^1 (\mu )$.   
  
\begin{coro}
\label{corergo}  
Let $(\theta (n) ; n\in \NN )$ be a sequence of complex numbers and   
$(u_n ; n\in \NN )$ be a strictly increasing sequence of integers. Suppose   
for some $\delta < 1$  we have $\H_1$ and for some integer $\gamma \geq \frac{1}{1-\delta   
}$  we have  
$$  
\sup_{N\geq 1} \frac{1}{N^{\gamma - 1}}     
\sum_{N^{\gamma} \leq k <(N+1)^{\gamma }} \left| \theta (k) \right| < \infty  .   
\leqno{( \H_3)}   
$$  
Then, the conclusions of Theorem \ref{thergo} hold.  
\end{coro}  
{\bf Proof.} It follows the lines of the proof of Theorem \ref{thergo}.  
\hfill $\Box$  
  
\medskip  
  
Condition $\H_3$ is useful in the case where the sequence $\theta$ is given by a random process (see \cite{DS}).  
  
Now under Condition $\H_2$ we give more information about the rate of convergence in Theorem \ref{thergo} when $f$ belongs to $L^2 (\mu )$. 

\begin{prop}
Let $\theta = (\theta (n) ; n\in \NN )$ be a bounded sequence of complex numbers and   
$u = (u_n ; n\in \NN )$ be a strictly increasing sequence of
 integers. Suppose Condition ${\H_2}$ holds.

Then, for all dynamical systems $(X,\B , \mu , T)$, all $f\in L^2 (\mu)$ and all $\beta > \rho $ we have 
$$  
\mu \left\{ x\in X ; \lim_{N\rightarrow +\infty } \frac{1}{N^{\beta}}   
\sum_{k=0}^{N-1} \theta (k) f (T^{u_k} x)  = 0 \right\} = 1 .  
$$  
\end{prop}
{\bf Proof.}
Let $f\in L^{2} (\mu )$ and $\beta > \rho $. From the Spectral Lemma (see \cite{Kr}) for all $m\leq n$ we have   
$$
\left| \left|   
\sum_{k=m}^{n} \theta (k) f \circ T^{u_k}   
\right| \right|_{2,\mu}  
=  
\left| \left|  
\sum_{k=m}^{n} \theta (k) e( \alpha u_k )  
\right| \right|_{2, \pTT,\mu_f ({\rm d}\alpha )}   
\leq   
D (n-m)^{\rho } \zeta (n)
\left| \left|  
f  
\right| \right|_{2,\mu} ,  
$$  
where $\mu_f$ is the spectral measure of $T$ at $f$. Then we proceed as
in the proof of Theorem \ref{thergo} with $\gamma > 1/(1-\rho )$.  Let
$M > 0$ and let $N$ be the unique integer such that $[N^{\gamma}] \leq M <   
[(N+1)^{\gamma}$]. Then  
$$  
\left|   
\frac{1}{M^{\beta}} \sum_{k=0}^{M-1} \theta (k) f (T^{u_k} x)   
\right|  
\leq  
\left|  
\frac{1}{N^{\beta \gamma}} \sum_{0\leq k < N^{\gamma }} \theta (k) f   
(T^{u_k} x)   
\right|  
$$
$$
+  
\left| \frac{1}{N^{\beta \gamma }} \sum_{N^{\gamma } \leq k < M-1}   
 \theta (k)  f (T^{u_k} x) \right| .  
$$
The first term tends to 0 and we have
$$
\left| \left|
 \frac{1}{N^{\beta \gamma }} \sum_{N^{\gamma } \leq k < M-1}   
\theta (k) f \circ T^{u_k}  
\right| \right|_{2,\mu}
\leq
D \frac{(M - N^{\gamma} )^{\rho }}{ N^{\beta \gamma } } \zeta (M)
\left| \left|  
f  
\right| \right|_{2,\mu}
$$
$$ 
\leq
D B^{\rho}N^{(\gamma - 1)\rho  - \beta\gamma} \zeta (M)
\left| \left|  
f  
\right| \right|_{2,\mu}
=
D B^{\rho}N^{\gamma(\rho - \beta) - \rho} \zeta (M)
\left| \left|  
f  
\right| \right|_{2,\mu} 
$$
which tends to 0 ($B$ is as in the proof of Theorem \ref{thergo}). This ends the proof.
\hfill $\Box$

\subsection{A precision in the case of the rotations}   
  
\label{cororot}

The goal of the following proposition is to investigate the uniform
convergence properties of the previous ergodic means. We look at the
particular case of the rotation dynamical systems. Let $A(\TT)$ be the
set of functions with summable Fourier coefficients, where $\TT$ is the
one-dimensional torus.

\begin{prop}
\label{rotation}  
Let $(\TT,\B,\mu , R_{\alpha })$ be a rotation dynamical system where   
$\alpha $ is an irrational number. Under the assumptions of Theorem   
\ref{thergo}, for all $f\in A (\TT)$ we have   
$$  
\frac{1}{N^{\beta}} \sum_{k=0}^{N-1} \theta (k) f (R_{\alpha}^{u_k} x)    
\longrightarrow_{N\rightarrow +\infty } 0    
$$  
uniformly in $x$ for all $\beta > \delta$.  
\end{prop} 

Let $f\in \A (\TT)$, i. e. $f (x) = \sum_{j\in \pZZ} \hat{f} (j) e(jx)$, $x\in \TT$,  with $\sum_{j\in \pZZ} |\hat{f}(j)| < \infty$. Let $N\in \NN$, we have  
$$  
\sum_{k=0}^{N-1} \theta (k) f \circ R_{\alpha}^{u_k} (x)  
=  
\sum_{j\in \pZZ} \hat{f} (j) e(jx) K_N (j\alpha),  
$$  
where $K_N (\lambda ) = \sum_{k=0}^{N-1} \theta (k) e(\lambda u_k
)$. Consequently Condition ${\cal H}_1$ and the fact that the map
$\alpha \mapsto j\alpha\ \mod \ 1 $ is onto imply that   
$$  
\sup_{x\in \pTT} \left|  \sum_{k=0}^{N-1} \theta (k) f \circ R_{\alpha}^{u_k} (x) \right|   
\leq    
CN^{\delta} \sum_{j\in \pZZ} \left|  \hat{f} (j) \right| ,  
$$  
which ends the proof.   
\hfill $\Box$  
  
\medskip  
  
Now we make a remark in the case we do not have a dynamical system structure, i.e., we are interested in the sequence $( \sum_{k=0}^{N-1} \theta (k) f (xu_k) ; N \in \NN) $. If $\theta $ satisfies Condition $\H_1$ and $f$ belongs to $\A (\TT)$ then the same computation as before leads to  
$$  
\sup_{x\in \pRR}   
\left|  
\sum_{k=0}^{N-1} \theta (k) f(xu_k)  
\right|  
\leq  
C N^{\delta } ||f||_{\A (\pTT)} .  
$$

\subsection{Proof of Theorem \ref{cvlaw}}

We first establish the following result.
  
\begin{theo}  
Let $H\in ]0,1[$. Let $(u_k ; k\in \NN)$ be a strictly increasing sequence such that there exist $\zeta \in ]0,1]$ and $\sigma > 1-H $ for which we have   
$$  
\left| \frac{u_N}{N} - \zeta \right| = O ( N^{-\sigma} ) . \leqno{(\H_4)}  
$$  
Then there exists a function $f\in L^2 (\mu)$ such that   
$$  
\frac{1}{N^{H}} \sum_{k=1}^{N} f\circ R_{\alpha}^{u_k} \Longrightarrow^{\hskip -12pt \L}_{\hskip +7pt } \N (0,1)   
.   
$$  
Moreover if $d$ is the Diophantine type of $\alpha$, then we can choose
 $f$ to be in ${\rm lip} (a)$ with $a<(1-H)/d$.   
\end{theo}  
{\bf Proof.}
From $\H_4$ we obtain the following estimation  
$$  
\sup_{\alpha \in \pRR}  
\left|  
\sum_{k=0}^{u_N} e(\alpha k) - \sum_{k=0}^{[\zeta N]} e(\alpha u_k)  
\right|  
= O (N^{1-\sigma }) .  
$$  
We set  
$$  
\Delta_N   
=  
\left|\left|  
\frac{1}{N^{H }} \sum_{k=0}^{u_N} f\circ R_{\alpha}^k - \frac{1}{N^{H}} \sum_{k=0}^{[\zeta N]} f\circ R_{\alpha}^{u_k}  
\right|\right|_{2,\mu} .  
$$   
From the Spectral Lemma we get   
$$  
\Delta_N  
\leq  
\frac{||f||_{2,\mu}}{N^{H}} \sup_{\alpha \in \pRR}  
\left|  
\sum_{k=0}^{u_N} e(\alpha k) - \sum_{k=0}^{[\zeta N]} e(\alpha u_k)  
\right|  
= O (N^{1-H -\sigma }) .  
$$  
But $1-H -\sigma $ is negative, hence $\Delta_N$ converges to
0. Now from a result of Lacey (Theorem 1.1 in \cite{La}) there exists a function $f\in L^2 (\mu)$ so that   
$$  
\frac{1}{N^{H}} \sum_{k=0}^{u_N} f   
\circ R_{\alpha}^{k}   \Longrightarrow^{\hskip -12pt \L}_{\hskip +7pt } \N (0,1)   
.    
$$
Moreover if $d$ is the Diophantine type of $\alpha$, then Theorem 1.1 in
\cite{La} allows us to choose $f$ to be in ${\rm lip} (a)$ with
$a<(1-H)/d$. 

We conclude the proof applying Slutsky Theorem. \hfill $\Box$

\medskip

We prove Theorem \ref{cvlaw} for $\theta^+$. No new arguments are needed to prove it for $\theta^-$. Let $\beta > \delta$. Let $N\in \NN$, we have  
  
$$  
\frac{1}{N^{\beta}} \sum_{k=0}^{N-1} \theta^{+} (k) f\circ  R_{\alpha}^{u_k}  
=  
\frac{1}{N^{\beta}} \frac{1}{2} \sum_{k=0}^{N-1} \theta (k) f\circ  R_{\alpha}^{u_k}  
+  
\frac{1}{N^{\beta}} \frac{1}{2} \sum_{k=0}^{N-1} f\circ  R_{\alpha}^{u_k} .  
$$  
  
From the Spectral Lemma and Condition $\H_1$ the first term  goes to
0 with respect to $|| . ||_{2,\mu}$, hence in probability. Theorem 1.1 in \cite{La} implies there exists $\tilde{f}$ in ${\rm lip} (a)$, with $a<(1-\beta)/d$, such that the second term converges in law to $\N (0,1)$. Consequently by Slutsky Theorem  
$$  
\frac{1}{N^{\beta}} \sum_{k=0}^{N-1} \theta^{+} (k) \tilde{f}   
\circ R_{\alpha}^{u_k}   \Longrightarrow^{\hskip -12pt \L}_{\hskip +7pt } \N (0,1)   
.    
$$  
This concludes the proof of Theorem \ref{cvlaw}.

\subsection{Discussion about Condition ${\cal H}_1$} 
 
Let us consider a strictly increasing sequence $u = (u_k ; k\in \NN )$ and a sequence $(\theta_k ; k\in \NN)$ fulfilling Condition ${\cal H}_1$, i.e., there exists $\delta <1$ such that $S(N,u) = O (N^{\delta})$. Let $(\beta_k ; k\in \NN)$ be an increasing sequence such that there exists $\gamma <1$ with $\beta_k = O(k^{\gamma})$. Then the sequences $\tilde{\theta} = (\theta_{k+\beta_k} ; k\in \NN)$ and $\tilde{u} = (u_{k+\beta_k} ; k\in \NN )$ satisfy Condition ${\cal H}_1$,  in fact for $\tilde {\delta} = \hbox {Max} (\delta \, , \gamma)$ we have

$$
S_N (\tilde{\theta} , \tilde{u}) := \sup_{\alpha \in
\pRR} \left|   \sum_{k=0}^{N-1} \tilde{\theta} (k) \exp (2i\pi \alpha \tilde{u}_k )
\right| = O (N^{\tilde{\delta}}) \ ,
$$ 

To show this it suffices to remark that
$$
\left \vert
\sum_{k=0}^{N+ \theta _N}\theta_ke(\alpha u_k)
- \sum_{k=0}^{N}\theta_{k+\beta_k}e(\alpha u_{k+\beta_k})
\right \vert \leq N^{\gamma}
\, .
$$

This remark allows us to construct deterministic sequences of weights
satisfying Condition ${\cal H}_1$ with $u_k = k + [\log (k+1)]$ for
example, where $[.]$ is the integer part map. In fact when $u= (k ; k\in
\NN )$, the sequence $\tilde{u}$ satisfies a condition of type $\H_4$.

A definition of the Thue-Morse sequence $\theta = (\theta_n ; n\in \NN)$ is the following. For all $n\in \NN$, let $r(n)$ be the sum modulo 2 of the digits of the expansion of $n$ in base 2, then $\theta$ can be defined by $\theta (n) = (-1)^{r(n)}$. 
We said in the introduction that $\delta (\theta , (k;k\in \NN))$ is
equal to $(\log 3)/ (\log 4)$ \cite{G}. Hence from what we said before
it comes that $\tilde{\delta}$ is also less than $(\log 3)/ (\log 4 )$, where $\tilde{\theta}_n = ( -1 )^{r(n+[\log (n+1)])}$ and $u_n = n+[\log (n+1)]$.

\section{$q$-multiplicative sequences}  
  
The goal of this subsection is to give some example of sequences
satisfying Condition $\H_1$. We mainly focus on $q$-multiplicative
sequences. We recall some known facts and results about these sequences and give a sufficient condition for a $q$-multiplicative sequence to fulfill $\H_1$.

\subsection{Definitions, notations and background}  
  
Let $q$ be an integer greater than or equal to 2. A {\it $q$-multiplicative   
sequence} $\theta = (\theta(n) ; n\in \NN )$ is a sequence of elements of $\UU$   
(the multiplicative group of complex numbers of modulus 1) such that for all integers $t \geq 1$ we have:  
  
\bigskip  
  
\centerline{$\theta(aq^t + b) =   
\theta (aq^t) \theta (b)$ for all $(a,b) \in \NN^2$ with $b<q^t$.}  
  
\bigskip  
  
We remark that necessarily $\theta (0) = 1$. The sequence $\theta$ is completely determined by the values of $\theta ( jq^k )$ where   
$(j,k)$ belongs to $\{ 0, \cdots , q-1 \} \times \NN$. Indeed if   
$  
n = \sum_{k\in \pNN} j_k q^k ,  
$   
where $j_k \in \{ 0, \cdots , q-1 \}$, for all $k\in \NN$, then   
$$  
\theta (n) = \prod_{k\in \pNN} \theta (j_k q^k).  
$$  
We will call {\it skeleton of $\theta$} the sequence $( (\theta (q^n),
\theta (2q^n), \cdots , \theta ((q-1)q^n) ) ; n\in \NN)$. We remark
easily that any sequence of $\UU$ is the skeleton of some
$q$-multiplicative sequence.

For all integers $N>0$ and all real numbers $x$ we set  
$$  
V_N (x) = \sum_{n=0}^{N-1} \theta (n) e(nx) \ \ ({\rm with} \ e(x) = e^{2i\pi   
x}).  
$$  
  
In \cite{LMM} the authors proved the following propositions.  
  
\begin{prop}  
\label{ufini}  
The following statements are equivalent.  
  
\medskip  
  
$\imath )$ $\theta (\NN)$ is finite,  
  
$\imath\imath )$ $\theta (\NN)$ contains an isolated point,  
  
$\imath\imath\imath )$ There exist $r\in \NN$ and $n_0$ such that for all   
$n\geq n_0$ and all $t\geq 0$, $\theta (tq^n)$ is an $r$-th root of unity.  
\end{prop}  
   
\begin{prop}  
If a non-periodic $q$-multiplicative sequence $\theta$ takes its values in a   
finite subset of $\UU$, then $\theta $ has empty spectrum, i.e., for all $x$ in $\RR$   
$$  
\lim_{n\rightarrow +\infty } \frac{1}{N} V_N (x) = 0.  
$$  
\end{prop}  
  
We set  $S_N (x) = V_{q^N} (x)$, we have $S_{N+1} = A_N (x) S_N (x)$, where 
$$
A_N (x) = \sum_{j<q} \theta (jq^N) e (jq^N x),
$$ 
and consequently   
$$  
S_N (x) = \prod_{n=0}^{N-1} A_n (x) .
$$

\subsection{A condition to fulfill Condition $\H_1$}
  
Let $\theta $ be a finitely valued (i.e., $\theta (\NN )$ is finite)
$q$-multiplicative sequence. For all $N\in \NN$ and all $0\leq j<q$, we
set $\theta (jq^N) = e (b_{N,j})$ with $0\leq b_{N,j}<1$. We remark that $b_{N,0} = 0$ for all $N\in \NN$. We set  
$$  
B_N =(1 ,\theta (q^N), \cdots , \theta ((q-1)q^N) ) \ {\rm and} \ E_N(x) =   
(1,e(q^N x), \cdots , e((q-1)q^Nx)).  
$$  
  
From Proposition \ref{ufini} there exist $r$ and $n_0$ such that   
for all $n\geq n_0$ and all $j\geq 0$ the complex number $\theta (jq^n)$ is an $r$-th root of unity. It comes that $b_{n,j}$ belongs to the set $\{ 0, 1/r , \cdots , (r-1)/r \}$. Hence the sets $\{ B_N ; N\in   
\NN    
\}$ and $\{ (b_{N,0} , b_{N,1} , \cdots , b_{N,q-1} ) ; N\in \NN \}$   
are finite.   
  
Let $U,V\in \CC^q$, we denote by $U.V$ the usual scalar product in $\CC^q$.   
We remark that $A_N (x) = B_N . E_N (x) = B_N . E_0 (q^N x)$ and of course $|A_n (x)|$ is less than $q$ for all $x\in \RR$.  
  
\medskip  
  
{\bf Remark.} $|B_N.E_0 (x)| = q$ if and only if for all   
$0\leq j \leq q-1$ we have  
$$  
b_{N,j} + j x \equiv 0,  
$$  
where $r\equiv s$ means $\{ r\} = \{ s\}$, $\{ . \}$ being the fractional part. Then it comes that $|B_N.E_0 (x)| = q$ if and only if  
$$  
x \equiv - b_{N,1} \hbox{ and } b_{N,j} \equiv j b_{N,1} \hbox{ for all } 0\leq j\leq q-1.  
$$

When the equation $|B_N.E_0 (x)| = q$ has a solution there is a
unique solution belonging to $[0,1[$ namely $x_N = 1 - b_{N,1} $ ($|B_N.E_0 (x_N)| = q$). Of course $x_N$ is not
defined for all $N$. We set $M = \{ n \in  \NN ; \sup_{x\in [0,1[} |B_n . E_n (x)| = q \}$ and   
$$
I=
\{
n\in M ; n+1\in M , \ \  qx_n \equiv x_{n+1}
\} .
$$
  
From the above relations we deduce that $I$ is the set of integers $N$ such that  
$$  
b_{N,j} \equiv j b_{N,1} \hbox{ and } b_{N+1,j} \equiv jq b_{N,1} \hbox{ for all } 0\leq j\leq q-1. \leqno{(R)}  
$$  
We set $I_N = I \cap [0, .. , N-1]$ and we say $\theta$ satisfies {\bf Condition (C)} if there exist $\alpha < 1$ such that we have   
$$  
\limsup_{N\rightarrow +\infty} \frac{\card I_N }{N}   \leq \alpha.  \leqno{(C)} 
$$  
The following proposition together with the relations (R) provide an easy way to construct $q$-mul\-ti\-plicative sequences fulfilling Condition $\H_1$.

\begin{theo}
\label{taux}  
Let $\theta$ be a finitely valued $q$-multiplicative sequence. If Condition (C) holds for $\theta$ then there exist $0<\delta<1$ and a constant $K$ such that
$$
\sup_{x\in \pRR} |V_N (x)| \leq K N^{\delta} \ \ \hbox{for all } N \in \NN.  
$$
\end{theo}  
  
{\bf Proof.} 
The set $\{ A_n A_{n+1} ; n\in \NN \}$ is finite. We claim that, for   
all $n\not \in I_N$, we have $\sup_{x\in [0,1]} |A_n(x)   
A_{n+1}(x)|^{1/2} = s_n < q$. Indeed if $\sup_{x\in [0,1]} |A_n(x)   
A_{n+1}(x)|^{1/2} =   
q$, for some $x$, then $|A_n(x)| = q$ and $|A_{n+1}(x)| = q$. Then it   
comes that $q^n x\equiv x_n$, $q^{n+1} x \equiv x_{n+1}$ and consequently   
$n\in I_N$.  
  
We set $s = \sup_{n\not \in I} \sup_{x\in [0,1]} |A_n(x)   
A_{n+1}(x)|^{1/2}$. We have $s<1$.  
  
By Condition (C) there exist $\alpha < 1$ and $N_0$ such that for all $N\geq N_0$ we   
have $\# I_N \leq \alpha N$. Let $N$ be such that $N-2 \geq N_0$. We   
have

$$  
| V_{q^N} (x) | = 
|S_N (x)| = \prod_{n=0}^{N-1} |A_n (x)| = |A_0 (x)|^{1/2}|A_{N-1}   
(x)|^{1/2}\prod_{n=1}^{N-2} |A_n (x)A_{n+1} (x)|^{1/2}    
$$  
$$  
\leq   
q^2 \prod_{n\in I_{N-2}} |A_n (x)A_{n+1} (x)|^{1/2} \prod_{n\not \in   
I_{N-2}} |A_n (x)A_{n+1} (x)|^{1/2}  
\leq   
q^2 q^{\# I_{N-2}} s^{N-2-\# I_{N-2}}  
$$  
$$  
=  
 q^2 \left( \frac{q}{s}\right)^{\# I_{N-2}}  s^{N-2}  
\leq  
 q^2 \left( \frac{q}{s}\right)^{\alpha (N-2)} s^{N-2} = q^2  (   
q^{\alpha} s^{1-\alpha } )^{N-2}.  
$$  
Hence there exist $\beta < 1$ and a constant $K$ such that for all   
$N\in \NN$ we have $\sup_{x\in [0,1]} |S_N (x)| \leq K q^{\beta N}$.

Let $a_n a_{n-1} \cdots a_0$, be the expansion of $N$ in
base $q$ with $a_n \not = 0$, $N = \sum_{i=0}^{n} a_i q^i$. We have  
$$  
|V_N (x)|   
=   
\left|   
\sum_{n=0}^{N-1} \theta (n) e(nx)   
\right|  
\leq   
$$  
$$  
a_0  + \left| \sum_{i=1}^{n} \sum_{j=0}^{a_i - 1} \sum_{k=0}^{q^i-1}   
\theta (a_0+\cdots a_{i-1}q^{i-1}+jq^i + k) e\left((a_0+\cdots   
a_{i-1}q^{i-1}+jq^i   
+ k)x\right) \right|  
$$  
$$  
\leq  
a_0  +  \sum_{i=1}^{n} \sum_{j=0}^{a_i - 1}   
\left|  
S_i (x)  
\right|  
\leq  
a_0  +  qK \sum_{i=1}^{n} q^{\beta i}   
\leq  
a_0  + qK q^{\beta (n+1)}  
\leq  
a_0  + Kq^{\beta + 1} N^{\beta}.  
$$  
  
This concludes the proof. \hfill $\Box$

\medskip

From this proof we have $\delta (\theta , (k;k\in \NN))
\leq \alpha + (1-\alpha )(\log s)/\log q$.

\medskip 

In the case where $q=2$ and $r=2$ then the $b_{n,j}$'s belong to the set $\{ 0,1/2\}$. Consequently Condition (C) becomes: there exist $\alpha <1$ and $N_0\in \NN$ such that for all $N\geq N_0$ we have
$$
\frac{1}{N} \card \{ 0\leq n \leq N-1 \ ; \  b_{n+1 , 1} \equiv 2b_{n,1}\}\leq \alpha  . 
$$ 
Then it is not difficult to see that Condition (C) holds if and only if there exist $\alpha <1$ and $N_0\in \NN$ such that for all $N\geq N_0$ we have
$$
\frac{1}{N}
\card 
\{
0\leq n \leq N-1 \ ; \  b_{n , 1} = 0
\} 
\leq \alpha .
$$
In this case the value $s$ used in the proof is less than or equal to $\frac{4}{3^{3/4}}$. It gives 

$$
\delta (\theta , (k;k\in \NN)) \leq 0,82 - 0,18	\alpha .
$$

Similar conditions can be given in the general case.
  
\subsection{Relations between Conditions (C) and $\H_2$}

In this subsection we prove that a finitely valued $q$-multiplicative
fulfilling Condition (C) satisfies Condition $\H_2$ (Proposition \ref{h3}). We prove this following the proof of Corollaire 1.11 in \cite{LMM}.

\begin{lemme}
Let $\theta$ be a $q$-multiplicative sequence. Then for all positive integers $N$, $p$ and $t$ we have for all $x\in \RR$
$$
\left|\sum_{n=0}^{N-1} \theta (n+p) e((n+p)x)\right| \leq 2 q^t + \frac{N}{q^t} \left| \sum_{n=0}^{q^t -1} \theta (n)e(nx) \right| .
$$
\end{lemme}
{\bf Proof.}
We set $a= [p/q^t]$ and $b=[(N+p)/q^t]$, where $[.]$ is the integer part map. If $N\leq 2q^t$, the inequality is clear. Otherwise we have $0\leq b-a\leq N/q^t$ and
$$
\left|\sum_{n=0}^{N-1} \theta (n+p) e((n+p)x)\right| \leq 2q^t + \left|\sum_{n=aq^t}^{bq^t - 1} \theta (n) e(nx)\right| 
$$
$$
= 2q^t + \left|\sum_{j=a}^{b - 1} \theta (jq^t)e(jq^t)  \right|\left|\sum_{k=0}^{q^t - 1} \theta (k) e(kx)\right|
\leq
2q^t + (b-a)\left|\sum_{k=0}^{q^t - 1} \theta (k) e(kx)\right|
$$
which ends the proof. \hfill $\Box$

\begin{prop}
\label{h3}
Let $\theta$ be a $q$-multiplicative sequence such that $\theta(\NN)$ is finite. If Condition (C) holds then there exist $0<\alpha<1$ and a constant $C$ such that
$$  
\sup_{x\in \pRR} \left| \sum_{n=0}^{N-1} \theta (n+p) e((n+p)x) \right| \leq C N^{\alpha} \  \hbox{ for all } N\in\NN \hbox{ and } p\in \NN.  
$$
\end{prop}  
{\bf Proof.}
We know there exist $0<\delta <1$ and a constant $C$ such that 
$$  
\sup_{x\in \pRR} \left| \sum_{n=0}^{N-1} \theta (n) e(nx) \right| \leq C N^{\delta} \ \ \hbox{ for all } N \in \NN.  
$$
Let $t$ be the unique integer such that 
$$
N^{\frac{1+\delta}{2}} \leq q^t < q N^{\frac{1+\delta}{2}} .
$$
Then the previous lemma gives
$$
\left| \sum_{n=0}^{N-1} \theta (n+p) e((n+p)x) \right|
\leq
2 qN^{\frac{1+\delta}{2}} + \frac{N}{N^{\frac{1+\delta}{2}}} Cq^{\delta} N^{\delta} = (2q+Cq^{\delta})N^{\frac{1+\delta}{2}} ,
$$
which ends the proof. 
\hfill $\Box$

\subsection{Some particular examples of $q$-multiplicative se\-quen\-ces}  

Before giving some examples, we need to recall some definitions about combinatorics on words.

An {\it alphabet} $A$ is a finite set of elements called {\it
letters}. A {\it word} on $A$ is an element of the free monoid generated
by $A$, denoted by $A^*$. Let $x = x_0x_1 \cdots x_{n-1}$ (with $x_i\in
A$, $0\leq i\leq n-1$) be a word, its {\it length} is $n$ and is denoted
by $|x|$. The {\it empty word} is denoted by $\mvide$, $|\epsilon| =
0$. The set of non empty words on  $A$ is denoted by $A^+$. The elements of $A^{\pNN}$ are called {\it sequences}. If $\x=\x_0\x_1\cdots$ is a sequence (with $\x_i\in A$, $i\in \NN$), and $I=[k,l]$ an interval of $\NN$ we set $\x_I = \x_k \x_{k+1}\cdots \x_{l}$ and we say that $\x_{I}$ is a {\it factor} of $\x$. If $k = 0$, we say that $\x_{I}$ is a {\it prefix} of $\x$. The {\it occurrences} in $\x$ of a word $u$ are the integers $i$ such that $\x_{[i,i + |u| - 1]}= u$. When $\x$ is a word, we use the same terminology with similar definitions.

The sequence $\x$ is {\it ultimately periodic} if there
exist a word $u$ and a non empty word $v$ such that $\x=uvvv\cdots $. Otherwise we say that $\x$ is {\it non-periodic}. It is {\it
periodic} (or {\it $|v|$-periodic}) if $u$ is the empty word.

The set $A$ is endowed with the discrete topology and $A^{\pNN}$ with the product topology. If $(u_n ; n\in \NN)$ is a sequence of words of $A^{*}$ such that $\lim_{n\rightarrow +\infty} |u_n| = +\infty $ then we say that $(u_n ; n\in \NN)$ converges to $u\in A^{\pNN}$ if and only if $(u_n^{\omega} ; n\in \NN)$ converges to $u\in A^{\pNN}$.   

\medskip 

{\bf Generalized Thue-Morse sequences}  
  
Let $r\geq 2$ be an integer and $R_r$ be the set of the $r$-th roots of   
 unity. We consider $R_r$ as an alphabet. Let $a_1\cdots a_n \in   
R_r^{*}$ and $b\in R_r$. We define $(a_1\cdots a_n ) * b$ to be the word $u$ of length $n$ defined by $u = (a_1.b) (a_2 . b) \cdots (a_n.b)$ where $x.y$ is the standard multiplication in $\CC$. The word $u$ belongs to $R_r^{*}$. In the same way we define   
$(a_1\cdots a_n )*(b_1\cdots b_m )$ to be the word $((a_1\cdots a_n   
)*b_1) \cdots ((a_1\cdots a_n )*b_m))$. It can be checked that this product is associative.

Let $(u_n ; n\in \NN)$ be a   
sequence of blocks of $R_r^{*}$, with $|u_n|\geq 2$, all beginning with the   
letter $1$, then the sequence of words $ (u_1 * u_2 * \cdots *  u_n ;   
n\in \NN ) $ converges to a sequence $\x\in R_r^{\pNN}$. We call it a   
{\it generalized Thue-Morse sequence}. These sequences were defined in \cite{Ke} for $r=2$. We will say it is of {\it constant length} whenever $|u_n| = |u_{n+1}|$ for all $n\in \NN$.  
  
\begin{prop}  
\label{gmorse}
Let $r\geq 2$ be an integer and $\x$ be a sequence of $R_r^{\pNN}$.   
Then, $\x$ is $q$-multiplicative if and only if $\x$ is a generalized   
Thue-Morse sequence of constant length $q$.  
\end{prop}   
{\bf Proof.} Let $\theta$ be a $q$-multiplicative sequence with skeleton $( s_n  ; n\geq 0) $, then it can be checked that the sequence $ (s_1 * s_2 * \cdots *  s_n ; n\in \NN ) $ converges to $\theta$ and conversely.
\hfill $\Box$  
  
\medskip  
  
The Thue-Morse sequence $\x$  was defined in \cite{Mo} to be the limit of
$(U_n ; n\in \NN)$ where $U_n \in \{ a,b \}^{*}$ is defined by $U_0 =
a$, $V_0 = b$, $U_{n+1} = U_n V_n$ and $V_{n+1} = V_n U_n$. Let $a=1$
and $b=-1$, then the sequence $( (1 \ (-1))^{*n} ; n\in \NN )$ converges
to the Thue-Morse sequence, where $(1 \ \ (-1))^{*n}$ is the $n$-th $*$-power
of the word $(1 \ (-1))$. We recall that in \cite{G} it is proved that
$\delta (\x , (k;k\in \NN))$ is equal to $(\log 3)/ (\log 4 )$.

Let $\theta$ be a $2$-multiplicative sequence on the alphabet $R_3 = \{ 1, j , j^2 \}$. Condition (C) is: there exist $\alpha <1$ and $N_0\in \NN$ such that for all $N\geq N_0$ we have
$$
\frac{1}{N}
\card 
\left\{
0\leq n \leq N-1 \ ; \  (b_{n , 1} , b_{n+1 , 1} ) \in \left\{ (0,0) , (1/3 , 2/3) , (2/3 , 1/3 ) \right\}
\right\} 
\leq \alpha .
$$

For example, let $\theta$ be the $2$-multiplicative sequence with the
periodic skeleton $(1,j),(1,j), \cdots $. From Proposition \ref{gmorse} we have  
$$  
u = (1\ j) * (1 \ j) * (1\ j) \cdots  
=  
(1 \ j \ j \ j^2) * (1 \ j) \cdots  
$$  
$$=  
(1 \ j \ j \ j^2 \ j \ j^2 \ j^2 \ 1) * (1 \ j) \cdots .  
$$  
We clearly have $\alpha = 0$. An elementary computation with Maple gives that 
 $\delta (\theta , (k;k\in \NN))$ is less than $0,93$.

\medskip

{\bf Substitutions and $q$-multiplicative sequences}  

A {\it substitution} on the alphabet $A$ is a map $\sigma : A \rightarrow A^{+}$. Using the extension by concatenation to words and sequences, $\sigma $ can be defined on the sets $A^{*}$ and $A^{\pNN}$. If for some letter $a$ the word $\sigma (a)$ begins with the letter $a$ and that $\lim_{n\rightarrow +\infty} |\sigma^n (a)| = +\infty$ then the sequence $( \sigma^n (a) ; n\in \NN)$ converges to a sequence $\x$ which satisfies $\sigma (\x) =\x$ : $\x$ is a {\it fixed point} of $\sigma$. We say $\sigma $ is of {\it constant length} $q$ whenever $|\sigma (b)| = q$ for all letters $b$ in $A$ (see \cite{Qu} for more details). 

\begin{prop}  
If the skeleton of a $q$-multiplicative sequence $\x$ is $n$-periodic then   
$\x$ is the fixed point of a substitution of constant length $q^n$.   
\end{prop}  
{\bf Proof.} The $*$-product being associative we can suppose that the skeleton of $\x$ is 1-periodic and equal to $((w_1 \cdots w_{q^n} ); n\in \NN)$. Let $\sigma : R_r \rightarrow R_r^{*} $ be the substitution of constant length $q^n$ defined   
for all $a\in R_r$ by $\sigma (a) = (a.w_1) \cdots (a.w_{q^n})$. Then $\x$   
is the fixed point of $\sigma$ starting with the letter 1.  
\hfill $\Box$  
  
\medskip  
  
In the previous example this gives $\sigma (1) = 1j$, $\sigma (j) = jj^2$   
and $\sigma (j^2) = j^2 1$.

\medskip  
  
\medskip  
  
{\bf Generalized Rudin-Shapiro sequences}  
  
Let $(u(n) ; n\in \NN)$ be the sequence where $u(n)$ is the number of   
blocks ``11'' in the binary expansion of the integer $n$. In \cite{AM}   
the authors proved for all $t\in \RR$ the sequence $v(t) = (e(tu(n)) ;   
n\in\NN))$ is such that $\delta (v(t)) < 1$ if and only if $t \not \in   
\ZZ $. And in particular that  $\delta (v(t), (k;k\in \NN)) = 1/2$ if $t \in \ZZ +   
1/2 $. In \cite{AL} is given a bunch of sequences for which $\delta
=1/2$ (see also \cite{MT}).
  
The sequence $v(1/2)$ is the well-known Rudin-Shapiro sequence which   
was the first example of a sequence with $\delta (v(t), (k;k\in \NN)) = 1/2$. It   
satisfies the following recurrence relation : $v_{2n}= v_{4n+1}= v_n$ and   
$v_{4n+3} = -v_{2n+1}$.

\section{Some precisions in the $q$-multiplicative case}

Let $\theta = (\theta (k) ; k\in \NN )$ be a non-periodic finitely valued $q$-multiplicative sequence. From Proposition \ref{ufini} there exists $n_0$ such that $\overline{\theta} = (\theta (k) ; k\geq n_0 )$ takes values in the set $S=\{ s_1 , \cdots , s_{|S|} \}$ contained in $R_r$, the set of the $r$-th roots of  unity for some $r\in \NN$. We define ${\rm Per} (\theta )$ to be the set of integer $a\in ]0,r[$ such that $(\theta^a (k) ; k\in \NN )$ is periodic. The sequence $\theta $ is said to be {\it irreducible} if ${\rm Per} (\theta ) = \emptyset$. 

\begin{prop}
Let $\theta = (\theta (k) ; k\in \NN )$ be a finitely valued $q$-multiplicative sequence. Let $\phi : S \rightarrow \CC$. For all dynamical systems $(X,\B,\mu, T)$ and all $f\in L^1 (\mu )$ the limit   
$$  
\lim_{N\rightarrow +\infty} \frac{1}{N} \sum_{k=0}^{N-1} \phi (\theta (k)) f\circ T^k  
$$  
exists and is equal to 0 if $\theta$ is irreducible.  
\end{prop}  
{\bf Proof.} We take the notations of the beginning of this subsection. In the sequel we suppose $\theta = \overline{\theta}$ but the same kind of proof holds when $\theta$ is not equal to $\overline{\theta}$.   
  
Let $(X,\B,\mu, T)$ be a dynamical system and $f\in L^1 (\mu )$. We take the notations of the beginning of this section.  
  
For all $s_j\in S$ we set $P_j (X) = \prod_{s\in R_r, s\not = s_j} (X - s) = a_{j,r-1}X^{r-1} + a_{j,r-2}X^{r-2} + \cdots +  a_{j,1}X + a_{j,0}$. We have  
$$  
\frac{1}{N}\sum_{k=0}^{N-1} \phi (\theta (k)) f\circ T^k (x)  
=  
\frac{1}{N}\sum_{j=1}^{|S|} \phi (s_j)\sum_{0\leq k \leq N-1 , \theta (k) = s_j} f\circ T^k (x)   
=  
$$  
$$  
=  
\frac{1}{N}\sum_{j=1}^{|S|} \phi (s_j) \sum_{k=0}^{N-1} \frac{\prod_{s\in R_r, s\not = s_j} (\theta (k) - s)}{\prod_{s\in R_r, s\not = s_j} (s_j - s)} f\circ T^k (x)  
$$  
$$  
=  
\frac{1}{N}\sum_{j=1}^{|S|} \frac{\phi (s_j)}{ \prod_{s\in R_r, s\not = s_j} (s_j - s) } \sum_{k=0}^{N-1}  P_j (\theta (k))    f\circ T^k (x)   
$$  
$$  
=  
\sum_{j=1}^{|S|}\sum_{l=0}^{r-1} \frac{\phi (s_j)a_{j,l}}{ \prod_{s\in R_r, s\not = s_j} (s_j - s) }\frac{1}{N} \sum_{k=0}^{N-1}  \theta^l (k)    f\circ T^k (x)   
$$  
If $l$ does not belong to ${\rm Per} (\theta )$ then $\theta^l$ is non-periodic and $q$-multiplicative. Consequently $\frac{1}{N} \sum_{k=0}^{N-1}  \theta^l (k)    f\circ T^k (x) $ converges to $0$.   
  
If $l$ belongs to ${\rm Per} (\theta )$ and set $\theta^l = yz_1 z_2 \cdots z_pz_1 z_2 \cdots z_p \cdots  $ where $y$ is a word on the alphabet $R_r $ and $z_i $ belongs to $ R_r$ for all  $1\leq i \leq p$. Then from Birkhoff's ergodic theorem  $\frac{1}{N} \sum_{k=0}^{N-1}  \theta^l (k)    f\circ T^k (x) $ converges to $(z_1 + \cdots + z_p) \int_X f d\mu$. \hfill $\Box$

\subsection{The case of the squares}  
  
We recall the van der Corput inequality (see \cite{KN}).  
  
\begin{lemme}  
Let $N\in \NN$ and $(\gamma (k) ; 0\leq k \leq N)$ be a finite sequence of a Hilbert space $\H$. For all $0\leq H\leq N-1$ we have  
$$  
\left| \left|  
\frac{1}{N}\sum_{k=0}^{N-1} \gamma (k)  
\right| \right|^2_{\H}  
\leq  
\frac{N+H}{N^2 (H+1)} \sum_{k=0}^{N-1} \left|\left| \gamma (k)\right|\right|_{\H}^2   
$$  
$$  
+   
\R {\rm e}   
\left(  
2\frac{N+H}{N^2 (H+1)^2}\sum_{h=1}^{H} (H+1-h)   
\sum_{k=0}^{N-h-1}   
< \gamma (k+h) , {\gamma (k)} >_{\H}  
\right)  
.  
$$    
\end{lemme}

\begin{theo}  
Let $\theta $ be a non-periodic finitely valued $q$-multiplicative
 sequence. Then, for all totally ergodic dynamical systems $(X,\B,\mu,T)$
 (i.e., $(X,\B,\mu,T^n)$ is ergodic for all $n\in \ZZ$), with $T$ invertible, and all $f\in L^2 (\mu )$ we have   
$$  
\lim_{N\rightarrow +\infty}  
\frac{1}{N} \sum_{k=1}^{N} \theta (k) f \circ T^{k^2} (x) = 0  
$$  
$\mu$-almost everywhere.  
\end{theo}  
{\bf Proof.} Let $f\in \H = L^2 (\mu )$ such that $\int_X f d\mu = 0$. We apply the van der Corput inequality to the sequence $(\gamma (k) ; k\in \NN )$ defined by $\gamma (k) = \theta (k) f\circ T^{k^2}$. For all $0\leq  H \leq N-1$ we obtain  
$$  
\frac{1}{N^2}  
\left|\left|  
\sum_{k=0}^{N-1} \gamma (k)  
\right|\right|_{\H}^2  
\leq  
\frac{2}{H+1} \left|\left| f \right|\right|_{\H}^2  
$$  
$$   
+  
\frac{2(N+H)}{N^2 (H+1)^2}\sum_{h=1}^{H} (H+1-h)   
\sum_{k=0}^{N-h-1}  
\R {\rm e}   
<  
\gamma (k+h), \gamma (k)   
>_{\H} \ \  
$$  
$$  
\leq  
\frac{2}{H+1} \left|\left| f \right|\right|_{H}^2  
$$  
$$  
+  
\R {\rm e}  
\left(  
\frac{2(N+H)}{N^2 (H+1)^2}\sum_{h=1}^{H} (H+1-h)   
\hskip -1.25pt  
\sum_{k=0}^{N-h-1}  
\hskip -1.25pt  
\theta (h+k) \overline{\theta (k)}  
<  
f\circ T^{(h+k)^2}, f\circ T^{k^2}   
>_{\H}  
\right)  
$$  
$$  
\leq \frac{2}{H+1} \left|\left| f \right|\right|_{\H}^2  
+   
\frac{4}{(H+1)^2}\sum_{h=1}^{H} (H+1-h) \left|  
<  
\frac{1}{N-h}  
\sum_{k=0}^{N-h-1}  
f\circ T^{ 2hk}, f \circ T^{-h^2}  
>_{\H}  
\right|   .
$$  
From the weak Ergodic Theorem and the facts that the system is totally
ergodic and $\int_X fd\mu = 0$, we have for all $H\in \NN$
$$  
\limsup_{N\rightarrow +\infty}   
\left|\left|  
\frac{1}{N}  
\sum_{k=0}^{N-1} \theta (k) f\circ T^{k^2}  
\right|\right|_{\H}^2  
\leq  
\frac{2}{H+1} \left|\left| f \right|\right|_{\H}^2 .  
$$  
It comes that $\lim_{N\rightarrow +\infty}   
\left|\left|  
\frac{1}{N}  
\sum_{k=0}^{N-1} \theta (k) f\circ T^{k^2}  
\right|\right|_{\mu , 2} = 0$.   
  
Now take $g\in L^2 (\mu )$. We apply what we just proved to $f = g -
\int_X g d\mu $. Because $\theta $ has an empty spectrum we obtain    
$$  
\lim_{N\rightarrow +\infty}   
\left|\left|  
\frac{1}{N}  
\sum_{k=0}^{N-1} \theta (k) g\circ T^{k^2}  
\right|\right|_{\mu , 2} = 0 .  
$$   
Moreover we know from \cite{Bo} that $(k^2 ; k\in \NN )$ is a good subsequence for the pointwise Ergodic Theorem in $L^2 (\mu)$ and furthermore from standard arguments we obtain the same conclusion for almost sure convergence.   
\hfill $\Box $  

\medskip

From the proof of the previous theorem and without using the fact that $(k^2
; k\in \NN )$ is a good subsequence for the pointwise Ergodic Theorem in
$L^2 (\mu)$, we obtain Proposition \ref{carre}.  
  
\bigskip
  
{\bf Proof of Proposition \ref{carre}.} When $\alpha $ is a rational
number a direct calculus leads to the result.

Now we suppose $\alpha$ is an irrational number. Let $\R_{\alpha}$ be the
rotation of angle $\alpha$ and $\mu$ be its Haar measure. It is easy to
check that it is totally ergodic.
From the previous proof it comes that for $\mu$-almost every $x\in [0,1[$ we have   
$$  
\lim_{N\rightarrow +\infty } \frac{1}{N} \sum_{k=0}^{N-1} \theta (k) e(\alpha k^2 + x ) = 0 .  
$$  
Hence $\lim_{N\rightarrow +\infty } \frac{1}{N} \sum_{k=0}^{N-1} \theta (k)   
e(\alpha k^2 ) = 0$ for all irrational $\alpha$. \hfill $\Box$

\bigskip

{\bf Acknowledgements.} We would like to thank Professor A.-H. Fan for
his valuable comments and suggestions. 

The main part of this work was done whence the first author was a member
of the Centro de Modelamiento Matem\'atico (UMR 2071 UCHILE-CNRS) in
Santiago de Chile and whence the second author was visiting the CMM. We
would like to acknowledge the CMM for its kind hospitality and for the
working environment it provided to us. 

We also thank the referee for his comments about Proposition \ref{carre}.

Fabien Durand, Centro de Modelamiento Matem\'atico, UMR 2071   
UCHILE-CNRS,  
Universidad de Chile, Casilla 170-3, Correo 3, Santiago, Chile, and 
 Laboratoire Ami\'enois de Math\'ematiques Fondamentales  et
Appliqu\'ees, CNRS-FRE 2270, Universit\'{e} de Picardie
Jules Verne, 33 rue Saint Leu, 80039 Amiens Cedex 01, France.

fdurand@u-picardie.fr
  
http://www.mathinfo.u-picardie.fr/fdurand/

\medskip  
  
Dominique Schneider, Laboratoire Ami\'enois de Math\'ematiques
Fondamentales  et Appliqu\'ees, CNRS-FRE 2270, Universit\'e de Picardie   
Jules Verne, 33, Rue Saint-Leu, 80039 Amiens Cedex 01, France.  

dominique.schneider@u-picardie.fr
  
http://www.mathinfo.u-picardie.fr/schneide/

\end{document}